\documentclass[11pt]{amsart}

\usepackage{}
\usepackage{mathrsfs}
\usepackage{amssymb}
\usepackage{amsmath}

\usepackage{color,latexsym,amsfonts}
\usepackage[active]{srcltx}
\usepackage[cp1252]{inputenc}
\usepackage{stmaryrd}

 \def\beqlb{\begin{eqnarray}}\def\eeqlb{\end{eqnarray}}
 \def\beqnn{\begin{eqnarray*}}\def\eeqnn{\end{eqnarray*}}

\newcommand{\bm}{{\mathbf m}}

\newcommand{\PF}{\mathbf {P}}
\newcommand{\EF}{\mathbf {E}}
\newcommand{\QZ}{{\mathbf Q_0}}

\newcommand{\ci}{{\mathcal I}}
\newcommand{\cj}{{\mathcal J}}
\newcommand{\ck}{{\mathcal K}}

\newcommand{\cn}{{\mathcal N}}

\newcommand{\Tau}{{\mathcal T}}
\newcommand{\ct}{{\mathcal T}}
\newcommand{ \cz}{\mathcal Z}

\newcommand{\E}{{\mathbb E}}
\newcommand{\N}{{\mathbb N}}
\renewcommand{\P}{{\mathbb P}}
\newcommand{\R}{{\mathbb R}}
\newcommand{\T}{{\mathbb T}}

\newcommand{\ind}{{\bf 1}}
\newcommand{\expp}[1]{\mathop {\mathrm{e}^{ #1}}}
\newcommand{\lb}{[\![}
\newcommand{\rb}{]\!]}

\newtheorem{thm}{Theorem}[section]
\newtheorem{lem}{Lemma}[section]

\newtheorem{cor}{Corollary}[section]
\newtheorem{defn}{Definition}[section]

\numberwithin{equation}{section}
\theoremstyle{remark}
\newtheorem{rem}{Remark}[section]

\def\pf{\noindent{\it Proof.~~}}
 \def\qed{\hfill$\Box$\medskip}

\begin{document}

\title{Time to MRCA for stationary CBI-processes}
\author{Hongwei Bi}

\address{School of Mathematical Sciences, Beijing Normal University, Beijing 100875,
P. R. China.}

\email{bihw@mail.bnu.edu.cn}

\thanks{This research  was supported by
  NSFC and 985 Program.}

\begin{abstract}
Motivated by sample path decomposition of
the stationary continuous state branching process with immigration, a
general population model is considered using the idea of immortal
individual. We compute the joint distribution of the random variables:
the time to the most recent common ancestor (MRCA), the size of the
current population and the size of the population just before MRCA.
We obtain the bottleneck effect as well.
The distribution of the number of the oldest families is also established.
The results generalize those in the recent paper by Chen and Delmas
\cite{cd10}.
\end{abstract}



\maketitle

\section{Introduction\label{s1}}

Continuous state branching processes (CB-processes) are non-negative
real-valued Markov processes first introduced by Jirina \cite{J58} to
model the evolution of large populations of small particles.
Continuous state branching processes with immigration (CBI-processes) are
generalizations of those describing the situation where immigrants may
come from outer sources, see e.g. Kawazu and Watanabe \cite{kw71}.
It is shown in Lamperti \cite{l67} that a CB-process can be obtained as the
scaling limit of a sequence of Galton-Watson processes; see also
\cite{Ali85, AlS82, Li06}. A genealogical tree is naturally associated
with the Galton-Watson process. This has given birth to the continuum
random tree theory first introduced by Aldous \cite{a91,a93} to code the
genealogy of the CB-process.  Duquesne and Le Gall \cite{dl02}
further developed the continuum L\'evy tree to give the complete
description of the genealogy of the CB-process in (sub)-critical case.
Kingman has initiated the study of the coalescent process
in 1982 in his famous papers \cite{k82a,k82b}.
Then coalescents with multiple collisions, also known as $\Lambda$-coalescents,
 were first introduced and studied independently by Pitman \cite{p99}
and by Sagitov \cite{s99}. Recently some authors have studied the
coalescent process associated with branching processes, see e.g. Lambert
\cite{l07} on coalescent time, Evans and Ralph \cite{er10} on the
dynamics of the time to the most recent common ancestor (MRCA),
Chen and Delmas \cite{cd10} on MRCA on some special stationary CBI-process,
and Berestycki et al.\cite{bbs11} on the coupling
between $\Lambda$-coalescents and branching processes.

This paper is motivated by Chen and Delmas \cite{cd10}.
The model considered here is a direct extension of \cite{cd10}.
We will use some notations and definitions in that paper and consider the general
CBI-process here instead. The fact that the CBI-process may have a non-trivial
stationary distribution makes it a more interesting model to be considered here
than the CB-process since for the CB-process either the population becomes
extinct or blows up with positive probability.
We consider a (sub)-critical CBI-process $Y = (Y_t, t \geq 0)$
with branching mechanism~$\psi$ given by \eqref{2.1} and immigration
mechanism~$F$ given by \eqref{2.7}. Our main interest is in presenting a further model of
random size varying population and exhibiting some interesting properties.
Afterwards we will give some properties of the coalescent tree.

We consider the stationary  CBI-process defined on the real line
$Z = (Z_t, t \in \R)$.
In order for the time to MRCA to be finite, we assume condition {\bf (A1)}:
$$
\int_1^{\infty} \frac{dz}{\psi(z)} < \infty.
$$
In order for $Z_t$ to be finite, we shall assume condition {\bf (A2)}:
$$
\int_0^\lambda \frac{F(z)}{\psi(z)}\, dz < \infty, \quad \mbox{for some}~
\lambda > 0.
$$

Using the look-down construction for the population with constant sizes,
we represent the process~$Z$ by means of the picture of an immortal individual
which gives birth to independent populations. We first give some
notations. For fixed time $t = 0$~(indeed we can choose any time by
stationarity), we denote by $A$ the time to the MRCA of the population
living at time 0, $Z^A = Z_{(-A)-}$ the size of the population just before
MRCA, and $Z^I$ the size of the population at time 0 which has been
generated by the immortal individual over the time interval $(-A,0)$ and
$Z^O = Z_0 - Z^I$ the size of the population at time 0 generated by the
immortal individual at time $-A$. We will see that conditionally on $A$,
the random variables $Z^A, Z^I$ and $Z^O $ are independent, and the joint
distribution of the random variables is also considered.
We also obtain the result that the size
of the population just before MRCA is stochastically smaller than that of
the population at the current time, that is the bottleneck effect.

Let $N^A + 1$ represent the number of individuals involved in the last
coalescent event of the genealogical tree.
We present the joint distribution of  $A, N^A$ and $Z_0$.
Using the measured rooted real tree formulation
of the genealogy of the stationary CBI-process developed in \cite{adh12},
we give the asymptotic for the number of ancestors.

We will give the transition probabilities of the MRCA age process $(A_t, t \in \R)$,
which has been studied by Evans and Ralph in \cite{er10}
for the CB-process conditioned on non-extinction.
We generalize it to the general case with the similar lines as their proof.
In the end we study the zero set of the CBI-process as well,
which is a stationary regenerative set.
Foucart and Bravo \cite{fb12} have studied the CBI case on the positive half line.
The stationary case is a bit different as the subordinator is not naturally
associated with the regenerative set; see \cite{t80} and \cite{ft88} for details.
For this situation see also \cite{ffs85}.

This paper is organized as follows. We first recall some well-known
results on the CB-process and CBI-process in Section 2. The family and clan
decomposition of the CBI-process are then introduced in Section 3. We
will give the condition for the existence of the stationary CBI-process,
determine the joint distribution of $A,Z^A,Z^I,Z^O$ and prove the
bottleneck effect in Section 4, that is $Z^A$\ is stochastically smaller than $Z_0$.
In Section 5 the distribution of the number of individuals involved in the last
coalescent event $N^A$ is computed. In the latter part of Section 6 we will
introduce the genealogy of CB-process using continuum random L\'evy
trees. Then  the asymptotic for the number of ancestors is given.
In Section 7 we give the transition probabilities
of the MRCA age process and the properties of the zero set.

\section{CB-process and CBI-process \label{s2}}

We recall some well-known results on CB-process and CBI-process derived
from Li \cite{l11, l12}. We consider a (sub)-critical branching mechanism
$\psi$:
\beqlb\label{2.1}
\psi(z) = b z + c z^2 + \int_0^{\infty} (\expp{-z u} - 1 + z u)\,m(du),\qquad z \geq 0,
 \eeqlb
where $b = \psi^{'}(0+) \geq 0, c \geq 0$ are constants and $(u \wedge u^2)\,m(du)$
is a finite measure on $(0, \infty)$. We will consider the non-trivial case, that is,
 assumption {\bf (A3)}:
$$
 \mbox{either}~c > 0~ \mbox{or}~ \int_{(0,1)} u\, m(du) = \infty.
$$
There exists a $\R_+$-valued strong Markov process $X = (X_t, t \geq 0)$
called continuous state branching process (CB-process) with branching mechanism $\psi$
whose distribution is characterized by its Laplace transform
 \beqlb
\EF_x[\expp{-\lambda X_t}] = \expp{-x\upsilon_t(\lambda)},
 \eeqlb
where $\EF_x$ means that $X_0 = x$ and the function $\upsilon_t(\lambda)$ is
the unique non-negative solution of the backward equation
 \begin{equation} \label{eq: bevtl}
\left\{  \begin{aligned}
&\frac{\partial}{\partial t} \upsilon_t(\lambda)
= -\psi(\upsilon_t(\lambda)),\quad t > 0,\lambda \geq 0,\\
&\upsilon_0(\lambda) = \lambda,\qquad\qquad\qquad  \lambda \geq 0,
 \end{aligned} \right.
 \end{equation}
 The CB-process has a canonical Feller realization.
 Let $\PF_x$ be the law of such a CB-process started at mass $x > 0$.
 Moreover, $X$ has no fixed discontinuities. The probability measure
 $Q_t(x, \cdot)$ is infinitely divisible and under condition {\bf (A3)},
 $\upsilon_t(\lambda)$ can be expressed canonically as
$$
\qquad \upsilon_t(\lambda)
= \int_0^{\infty} (1- \expp{-\lambda u})\,l_t(du),\qquad  t > 0, \lambda \geq 0,
$$
where $u\, l_t(du)$ is a finite measure on $(0, \infty)$; see Theorem 3.10 in
\cite{l11}. The Markov property of $X$ implies that for any $\lambda, s, t \geq 0$,
 \begin{equation}\label{cumu-up}
\upsilon_{t + s}(\lambda) = \upsilon_t(\upsilon_s(\lambda)).
 \end{equation}
We also have the forward differential equation
 \beqlb \label{forwarde}
\frac{\partial}{\partial t}\upsilon_t(\lambda) =
-\psi(\lambda) \frac{\partial}{\partial \lambda} \upsilon_t(\lambda).
 \eeqlb
Let $\zeta = \inf\{s \geq 0, X_s = 0\}$ be the extinction time of $X$ and
$c(t) = \lim_{\lambda \to \infty} \upsilon_t(\lambda)$. Under {\bf (A1)},
$c(t) > 0$ is finite. We have by \eqref{cumu-up} that
\begin{equation}\label{cumu-ct}
\upsilon_s (c(t)) = c(t + s).
\end{equation}

We consider an immigration mechanism $F$:
 \beqlb\label{2.7}
F(z) = \beta z + \int_0^{\infty} (1 - \expp{-z u})\,n(du),\qquad z \geq 0,
 \eeqlb
where $\beta \geq 0$ is a constant and $(1 \wedge u)\,n(du)$ is a finite
measure on $(0, \infty)$. Then there exists a strong Markov process
$Y = (Y_t, t \geq 0)$  called continuous-state
branching process with immigration (CBI-process) with branching mechanism
$\psi$ and immigration mechanism $F$ defined on $\R_+$ with Laplace transform
given by
\beqlb
\E_x[\expp{-\lambda Y_t}]
= \expp{-x \upsilon_t(\lambda) - \int_0^t F(\upsilon_s(\lambda))\,ds},
 \eeqlb
where $\E_x$ means that $Y_0 = x$.
We also denote $\P$ the corresponding probability measure.

\section{Sample path decomposition\label{s3}}

In this section we will recall some results from Li \cite{l11,l12}. We
will give the clan and family
decomposition of the CBI-process.

The process $X$ is infinitely divisible. It is well known that there exists
a canonical measure (we also call it excursion law) $\QZ$ on the space $D$
of C\`adl\`ag  functions on $[0, \infty)$
with Skorokhod topology.  Notice that $\QZ(\{X, X_{0+} \neq 0\}) = 0.$

We can give a reconstruction of the sample paths of the CB-process by
means of the excursion law. Let $x \geq 0$ and let
$N(dX) = \sum_{i \in I} \delta_{X^i}(dX)$ be a
Poisson random measure on $D$ with intensity $x \QZ(dX)$.
We define the process $(X_t^{'}, t\geq 0)$ by
 \beqlb\label{eq:dspem}
\left\{ \begin{aligned}
X_0^{'} &= x,\\
X_t^{'} &= \int_{D}X_t\,N(dX),\qquad t > 0.
 \end{aligned}\right.
 \eeqlb
Then $X^{'}$ is a realization of the CB-process $X$. We will not
distinguish $X^{'}$ from $X$. For the proof see Li \cite[Theorem 8.24]{l11}
 or \cite[Theorem 2.4.2]{l12}. As one can see \eqref{eq:dspem} is
equivalent to the well-known decomposition as follows:
If $N(dx, dX) = \sum_{i \in  I} \delta_{(x_i, X^i)} (dx, dX)$
is a Poisson point measure on $\R_+ \times D$
with intensity $\ind_{[0, \infty)}(x) dx \QZ (dX)$,
then $\sum_{i \in I} \ind_{\{x_i \leq x\}} X^i$  is distributed as $X$
under $\PF_x$. Further we have for $\lambda \geq 0$,
 \beqnn
\QZ(1 - \expp{-\lambda X_t}) = \lim_{x \to 0}\frac{1}{x}
\EF_x[1 - \expp{-\lambda X_t}] = \upsilon_t(\lambda),
 \eeqnn
and
 \beqnn
c(t) = \QZ(\zeta > t) = \QZ(X_t > 0).
 \eeqnn
We will put $X_t = 0$ for $t < 0$.

Now we will introduce the family decomposition of the CBI-process. We
consider the following Poisson point measures.
 \begin{enumerate}
\item Let $N_0(dr, dt) = \sum_{i \in I} \delta_{(r_i, t_i)}(dr, dt)$ be a
    Poisson point measure on $(0, \infty) \times \R$ with intensity $n(dr) dt$.
\item Conditionally on $N_0$, let $(N_{1,i}, i \in I)$ be independent Poisson point
    measures with intensity $r_i \delta_{t_i}(dt) \QZ(dX)$, where
    $N_{1, i}(dt, dX) = \sum_{j \in \cj_{1, i}} \delta_{(t_j, X^j)}(dt, dX)$.
    Note that for all $j \in \cj_{1, i}$, we have $t_j = t_i$.
    We set $\cj_1 = \bigcup_{i \in  I} \cj_{1, i}$,
    and $N_1(dt, dX) = \sum_{j \in \cj_1} \delta_{(t_j, X^j)}(dt, dX)$.
\item Let $N_2(dt, dX) = \sum_{j \in \cj_2} \delta_{(t_j, X^j)}(dt,dX)$ be a Poisson point
    measure with intensity $\beta dt$ $\QZ(dX)$  independent of $N_0, N_1$.
 \end{enumerate}

We set $\cj = \cj_1 \bigcup \cj_2$. We shall call $X^j$ a family and
$t_j$ its birth place for $j \in \cj$. We will consider the process
$(Y^{'}_t, t \geq 0)$ and its stationary version $(Z_t, t \in \R)$.
They are usually called family decomposition of the CBI-process defined as follows:
\beqlb \label{eq:familyd}
Y_t^{'} = \sum_{j \in \cj, t_j > 0} X_{t - t_j}^j,\quad
Z_t = \sum_{j \in \cj} X_{t - t_j}^j.
 \eeqlb

Putting this another way we can deduce that it corresponds to a special
immigration process shown in Corollary 3.4.2 in \cite{l12}.

For $i\in I$ denote $X^i = \sum_{j \in \cj_{1, i}} X^j$ and
$\ci = I \bigcup \cj_2$. The random measure
 \beqlb
N_3(dt, dX) = \sum_{i \in \ci} \delta_{(t_i, X^i)}(dt,dX)
 \eeqlb
is a Poisson point measure with intensity $dt \mu(dX)$,
where $\mu$ is given by
\beqlb
\mu(dX) = \beta \QZ(dX)
+ \int_{(0, \infty)} n(dx) \PF_x(dX).
\eeqlb
It corresponds to the entrance law $(H_t, t > 0)$  in Li \cite{l12} ~given by
 \beqnn
H_t = \beta l_t + \int_0^{\infty} n(dx)\, Q_t(x, \cdot).
 \eeqnn
We shall call $X^i$ with $i \in \ci$ a clan and $t_i$ its birth place.
For $j \in \cj_2$, $X^j$ is a clan and a family. Then we have
 \beqlb \label{eq:cland}
Y_t^{'} = \sum_{i \in \ci, t_i > 0} X_{t - t_i}^i, \quad
Z_t = \sum_{i \in \ci} X_{t - t_i}^i.
 \eeqlb
$Y^{'}$ with this representation is just the sample path decomposition of
the CBI-process. We shall call this the clan decomposition of $Y$.
$Y^{'}$ is a version of $Y$ and $Z$ is the stationary version of $Y$.
Usually the family decomposition is more precise than the clan decomposition.

We give an interpretation of~$Z$ in population terms. At time t, $Z_t$
corresponds to the size of the population generated by an immortal
individual giving birth at rate $\beta$ with sizes evolving independently
as $X$ under $\QZ$ and at rate 1 with intensity $n(dx)$
with initial size $x$ which evolve independently as $X$ under
$\PF_x$.
We first give a lemma on the family representation.
 \begin{lem}\label{l3.1}
Let $f$ be a non-negative measurable function. We have
 \beqlb
\E\left[\expp{-\sum_{j \in \cj} f(t_j, X^j)}\right] =
\expp{-\int_{\R} F \left(\QZ(1 - \expp{-f(t, X)})\right)\,dt}.
 \eeqlb
 \end{lem}
\pf Due to the independence of the Poisson random measures and
the exponential formula, we have
 \begin{align*}
\E\left[\expp{-\sum_{j \in \cj} f(t_j, X^j)}\right]
& = \E\left[\expp{-\sum_{j \in \cj_1} f(t_j,X^j)}\right]
\E\left[\expp{-\sum_{j \in \cj_2} f(t_j, X^j)}\right]\\
& = \E\left[\expp{-\sum_{i \in I} \sum_{j \in \cj_{1, i}}f(t_j, X^j)}\right]
\expp{-\beta \int_{\R} \QZ(1 - \expp{-f(t, X)})\,dt}\\
& = \E\left[\expp{-\sum_{i \in I} r_i \QZ(1 - \expp{-f(t_i, X)})}\right]
\expp{-\beta \int_{\R} \QZ(1 - \expp{-f(t, X)})\,dt}\\
& = \expp{-\int_{\R} F(\QZ(1 - \expp{-f(t, X)}))\,dt}.
 \end{align*}\qed

The necessary and sufficient condition for which the CBI-process has a stationary version is that
{\bf (A2)} holds, see Theorem 3.20 in \cite{l11}.
If {\bf (A2)} holds, then $X_t$ converges in distribution to $X_{\infty}$ as $t \to \infty$,
with the distribution of $X_{\infty}$ characterized by its Laplace transform
\beqlb \label{leta}
\E[\expp{- \lambda X_{\infty}}] = \expp{-\int_0^{\infty}F(\upsilon_s(\lambda))\,ds}.
 \eeqlb
Then with {\bf (A2)} in force, $Z$ defined by \eqref{eq:familyd}
and \eqref{eq:cland} is the stationary version of $Y$.
 \begin{cor} \label{c3.1}
Assume that {\bf (A2)} holds. We have for $\lambda > 0, t \in \R,$
 \beqlb
\E[Z_t \exp(-\lambda Z_t)]=
\frac{F(\lambda)}{\psi(\lambda)} \E[\expp{-\lambda Z_t}].
 \eeqlb
In particular, we have
 \beqlb
\E[Z_t] = \frac{F^{'}(0)}{\psi^{'}(0)} .
 \eeqlb
 \end{cor}
\pf We can see from \eqref{leta} that
 \beqnn
\E[Z_t \exp(-\lambda Z_t)] = \E[\expp{-\lambda Z_t}]
\partial_{\lambda} \int_0^{\infty} F(\upsilon_s(\lambda))\, ds.
 \eeqnn
Using the forward equation \eqref{forwarde} we can deduce that
 \begin{eqnarray*}
 \begin{aligned}
\partial_{\lambda} \int_0^{\infty} F(\upsilon_s(\lambda))\, ds
&= \int_0^{\infty} \partial_x F(x) |_{x = \upsilon_s(\lambda)}
\partial_{\lambda} \upsilon_s(\lambda)\, ds\\
&= -\frac{1}{\psi(\lambda)} \int_0^{\infty}
F^{'}(\upsilon_s(\lambda)) \partial_s\upsilon_s(\lambda)\, ds
= \frac{F(\lambda)}{\psi(\lambda)}.
 \end{aligned}
 \end{eqnarray*}
The second part is obvious. \qed

In the following we will always suppose that {\bf (A1)}, {\bf (A2)} and
{\bf (A3)} are in force.

\section{Time to MRCA and the population sizes\label{s4}}

With the decomposition procedure in force, we will follow the steps of
Chen and Delmas \cite{cd10}. We consider the coalescence of the genealogy
at a fixed time $t_0$. We may as well assume that $t_0=0$ because of
stationarity. There are infinitely many number of clans contributing to the
population at time 0. We can further prove that there are only finite number of
clans born before time $a$ and still alive at time 0. Only one oldest
clan is expected to be still alive at time 0.

First we will give the notations using the decomposition. $-A$ is the birth time of
the unique oldest clan at time 0 ($A$ is also the time to the most recent
common ancestor (TMRCA) of the population at time 0) given by
$A:= -\inf\{t_i \leq 0, X_{-t_i}^i > 0, i \in \ci\}$;
$Z^O$ is the population size of this clan at time 0, i.e.
$Z^O:= X_{-t_i}^i$, if $A = -t_i$;
The size of all the clans alive at time 0 with birth time in $(-A,0)$ is
given by $Z^I:= Z_0 - Z^O$,
and the size of the population just before the MRCA is given by
$Z^A:= Z_{(-A)-} = \sum_{i\in\ci} X^i_{-A - t_i}1_{\{t_i < -A\}}$.

 \begin{thm}\label{t4.1}
Let $f:\R \rightarrow \R_+$ be a measurable function.
For $\lambda, \gamma, \eta \geq 0$, we have
 \begin{eqnarray*}
 \begin{aligned}
&\E[(\expp{-\lambda Z^A - \gamma Z^I - \eta Z^O} f(A))]\\
&=\int_0^{\infty} dt\,f(t) \left(F(c(t)) - F(\upsilon_t(\eta))\right)
\exp\left(-\int_0^t F(\upsilon_s(\gamma))\,ds
- \int_0^{\infty} F(\upsilon_s(\lambda + c(t))) ds\right).
 \end{aligned}
 \end{eqnarray*}
 \end{thm}
\pf  We have
 \begin{eqnarray*}
 \begin{aligned}
&\E[\expp{-\lambda Z^A - \gamma Z^I - \eta Z^O} f(A)]\\
&\qquad = \E\bigg[\sum_{j \in \ci} \exp\bigg(-\lambda \sum_{i \in \ci, t_i < t_j}
X_{t_j - t_i}^i - \gamma \sum_{i \in \ci, t_i > t_j} X_{-t_i}^i - \eta X_{-t_j}^j\bigg) \\
&\qquad\qquad\qquad\qquad\qquad\qquad\times f(-t_j) \ind_{\{X_{-t_j}^j > 0,
\sum_{i \in \ci, t_i < t_j} \ind_{\{X_{-t_i}^i > 0\}} = 0\}}\bigg]\\
&\qquad = \int_0^{\infty} dt\, f(t) \mu(\expp{-\eta X_t} \ind_{\{X_t > 0\}})
\E\bigg[\exp\bigg(-\gamma \sum_{i \in \ci, t_i > -t} X_{-t_i}^i\bigg)\bigg]\\
&\qquad\qquad\qquad \times \lim_{K \to \infty}
\E\bigg[\exp\bigg(-\lambda \sum_{i \in \ci,  t_i < -t}
\big(X_{-t - t_i}^i + K \ind_{\{X_{-t_i}^i > 0\}}\big)\bigg)\bigg],
 \end{aligned}
 \end{eqnarray*}
where the first equality is based on the values of $A$ and the second one
holds since Poisson point measures over disjoint sets are independent. We
will calculate the terms respectively.

First we have
 \beqnn
\mu(\expp{-\eta X_t} \ind_{\{X_t > 0\}})
= \mu \big(\ind_{\{X_t > 0\}} - (1- \expp{-\eta X_t}) \big).
 \eeqnn
Using the expression of $\mu$, we have
 \begin{eqnarray*}
 \begin{aligned}
\mu(\ind_{\{X_t > 0\}})
&= \beta \QZ(X_t > 0) + \int_0^{\infty} n(dx)\,\PF_x(X_t > 0)\\
&= \beta c(t)+ \int_0^{\infty} n(dx)\,(1- \PF_x(X_t = 0))\\
&= \beta c(t)+ \int_0^{\infty} n(dx)\,(1- \expp{-x c(t)}) = F(c(t)),
 \end{aligned}
 \end{eqnarray*}
and
 \begin{eqnarray*}
 \begin{aligned}
\mu(1-\expp{-\eta X_t}) &= \beta \QZ(1 - \expp{-\eta X_t})
+ \int_0^{\infty} n(dx)\,\PF_x(1- \expp{-\eta X_t})\\
&= \beta \upsilon_t(\eta) + \int_0^{\infty} n(dx)\,(1- \expp{-x\upsilon_t(\eta)}) = F(\upsilon_t(\eta)).
 \end{aligned}
 \end{eqnarray*}

Second, using Lemma \ref{l3.1} we get
 \beqnn
\E\left[\exp\bigg(-\gamma\sum_{i\in\ci,t_i>-t} X_{-t_i}^i\bigg)\right]
=\exp\bigg(-\int_0^t F\big(\upsilon_s(\gamma)\big)\,ds\bigg).
 \eeqnn

Finally we see that
 \begin{eqnarray*}
 \begin{aligned}
&\lim_{K \to \infty} \E\left[\exp\bigg(-\lambda \sum_{i \in \ci,  t_i < -t}
(X_{-t- t_i}^i + K \ind_{\{X_{-t_i}^i > 0\}})\bigg)\right]\\
&= \exp\bigg(-\int ds\, \ind_{\{s > 0\}}
\mu\big(1 - \expp{-\lambda X_s} \ind_{\{X_{s + t} = 0\}}\big)\bigg)\\
&= \exp\bigg(-\int ds\,\ind_{\{s > 0\}}
\mu\big(1 - \expp{-\lambda X_s} \PF_{X_s}(X_t = 0)\big)\bigg)\\
&= \exp\bigg(-\int ds\,\ind_{\{s > 0\}}
\mu\big(1 - \expp{-(\lambda + c(t)) X_s}\big)\bigg)\\
&= \exp\bigg(-\int_0^{\infty} ds\,F(\upsilon_s\big(\lambda + c(t))\big)\bigg),
 \end{aligned}
 \end{eqnarray*}
where we use the exponential formula for the Poisson point measure in the
first equality and the Markov property of $X$ in the second one.

Putting all the calculations together we obtain the result.\qed

It is then straightforward to derive the distribution of the TMRCA $A$.
 \begin{cor}\label{c4.1}
The distribution of  $A$ is given by
 \beqnn
\P(A\leq t) = \exp\bigg(-\int_t^{\infty}
F(c(s))\,ds\bigg) = \E[\expp{-c(t) Z_0}],
 \eeqnn
and $A$ has density $f_A$ with respect to the Lebesgue measure given by
 \beqnn
f_A(t) = \ind_{\{t > 0\}} F(c(t)) \exp\bigg(-\int_t^{\infty} F\big(c(s)\big)\,ds\bigg).
 \eeqnn
 \end{cor}
\pf Using Theorem \ref{t4.1} and \eqref{cumu-ct},
we see that the first equality holds easily.
The second one is immediate. \qed

We see that in this general case the expression of the distribution of $A$
is invariant compared with \cite{cd10}.
The next result is also a direct consequence of Theorem \ref{t4.1}.

 \begin{cor}\label{c4.2}
Conditionally on $A$, the random variables $Z^I, \ Z^A, \ Z^O$ are independent.
 \end{cor}

We also derive from Theorem \ref{t4.1} the distribution and the mean of the population size just before MRCA.
As can be seen from below that the expression for the Laplace transform is the same as that of \cite{cd10}.
 \begin{cor}\label{c4.5}
Let $t > 0$. Then
 \beqlb \label{eq:zaa}
\E\big[\expp{-\lambda Z^A}| A = t \big]=
\frac{\E[\expp{-(\lambda + c(t)) Z_0}]}{\E[\expp{-c(t) Z_0}]} \quad
and \quad \E[Z^A | A = t] = \frac{F(c(t))}{\psi(c(t))}.
 \eeqlb
 \end{cor}
\pf This is a direct consequence of Theorems \ref{t4.1} and Corollary \ref{c4.1}.\qed

We can further deduce that conditionally on $\{A = t\}$, the distribution of $Z^A$
converges to the distribution of $Z_0$ as $t \to \infty$.

Another application of Theorem \ref{t4.1}, we call the bottleneck effect,
 is that the size of the population just before MRCA is stochastically
smaller than that of the current population.
Note that this inequality does not hold in
the almost surely sense in general. The proof is the same as that of \cite{cd10}.

 \begin{cor}\label{cbneck}
For all $z \geq 0$ and $t \geq 0$, we have
$\P(Z^A \leq z| A = t) \geq \P(Z_0 \leq z)$.
Hence the population size $Z^A$ is stochastically smaller than $Z_0$,
that is  $\P(Z^A \leq z) \geq \P(Z_0 \leq z)$, for all $z \geq 0$.
In particular we have $\E[Z^A | A] \leq \E[Z_0]$.
 \end{cor}

 \begin{rem}\label{r1}
Instead of considering the size of the population just before MRCA,
we consider the size  at MRCA, $Z^A_+$, which is given as
$Z^A_+ = Z^A + \sum_{i \in I} X_0^i \ind_{\{t_i = -A\}}$. We don't take
into account the contribution of~$i \in \cj_2$ since for those we have $X_0^i = 0$.
Similar calculations as those of Theorem \ref{t4.1} show that for $\lambda, t > 0$,
 \beqnn
\E[\expp{-\lambda Z^A_+} | A = t] = \E[\expp{-\lambda Z^A} | A = t]
\frac{F(\lambda + c(t)) - F(\lambda)}{F(c(t))}.
 \eeqnn
If $F^{'}(0) = \infty$, then $\lim_{t\to\infty} \E[\expp{-\lambda Z^A_+} | A = t]= 0$,
 which means that conditionally on $\{A = t\}$, $Z_+^A$ is likely to be very large,
 as $t \to \infty$. We can interpret it as this:
  a clan is born at time $-t$ and it survives up to time 0,
  if $t$ is large enough, it is likely to have a large initial size.
  Therefore, $Z_+^A$ is not stochastically smaller than $Z_0$ in general.
 \end{rem}

\section{The number of oldest families\label{s5}}

In this section we will consider the number of families in the oldest clan alive at time $0$.
 It is equivalent as that of individuals involved in the last coalescent event
 of the genealogical tree. We will use the family representation in this section.

 \begin{def}\label{d5.1}
The number of oldest families alive at time 0 (excluding the immortal individual) is defined as:
 \beqlb
N^A = \sum_{j \in \cj} \ind_{\{A = -t_j,~X_{-t_j}^j > 0\}}
= \sum_{j \in \cj}\ind_{\{A = -t_j,~\zeta_j > -t_j\}}.
 \eeqlb
Obviously $N^A \geq 1$. In particular when $\beta > 0$
and the measure $n \equiv 0$, we have $N^A = 1$.
 \end{def}

The following theorem gives the joint distribution of $A, N^A$ and $Z_0$.
 \begin{thm} \label{t5.1}
Let $0 \leq a \leq 1$. For any non-negative measurable function $f$,  we have
 \begin{eqnarray*}
 \begin{aligned}
\E\big[a^{N^A} \expp{-\lambda Z_0} f(A)\big]
&= \int_0^{\infty} ds\, f(s)
\exp\bigg(-\int_0^s  F(\upsilon_r(\lambda)) dr - \int_s^{\infty} F(c(r))\, dr\bigg)\\
&\qquad\qquad\qquad\times\bigg(F(c(s)) - F\big((1 - a) c(s) + a \upsilon_s(\lambda)\big)\bigg).
 \end{aligned}
 \end{eqnarray*}
 \end{thm}
\pf For $i\in\ci$, set
 \beqnn
J_i^{*}=\left\{
 \begin{array}{cl}
\cj_{1,i}, & \quad\textrm{if}~i\in I,\\
\{i\}, & \quad\textrm{if}~i\in\cj_2.\\
 \end{array}\right.
 \eeqnn
For $f$ non-negative measurable, we have
 \begin{align*}
&\E\big[a^{N^A} \expp{-\lambda Z_0} f(A) \big]\\
&= \E\bigg[\expp{-\lambda \sum_{k \in \ci,  t_k < 0 } X_{-t_k}^k}
\sum_{i \in \ci} a^{\sum_{j \in J_i^{*}} \ind_{\{\zeta_j > -t_i\}}} f(-t_i)
\ind_{\{X_{-t_i}^i > 0\}} \ind_{\big\{\sum_{l \in \ci, t_l < t_i}
\ind_{\{X_{-t_l}^l > 0\}} = 0 \big\}}\bigg\}\\
&= \int_0^{\infty}ds\, f(s) \E\bigg[\expp{-\lambda \sum_{k \in \ci, t_k < 0}
X_{-t_k}^k \ind_{\{t_k > -s\}}}\bigg]
\P\bigg(\sum_{k \in \ci} \ind_{\{t_k < -s,X_{-t_k}^k > 0\}} = 0\bigg)\\
&\quad \times \bigg(\beta \QZ[a \expp{-\lambda X_s} \ind_{\{X_s > 0\}}]
+ \int_0^{\infty} n(dx) \E_x \bigg[a^{\sum_{j \in J_3}
\ind_{\{X_s^j > 0\}}} \expp{-\lambda \sum_{j \in J_3} X_s^j}
\ind_{\{\sum_{j \in J_3} X_s^j > 0\}}\bigg]\bigg),
 \end{align*}
where the first equality is based on the decomposition of $A$,
the second on splitting $t_k$ into three parts: $t_k>s,t_k<s$ and $t_k=s$,
and $\sum_{j\in J_3}\delta_{X^j}(dX)$ is a Poisson point
measure with intensity $x\QZ(dX)$ under $\P$.
We will calculate the terms separately.

By Lemma \ref{l3.1}, we have
 $$
\E\bigg[\exp\bigg(-\lambda \sum_{k \in \ci, t_k < 0} X_{-t_k}^k \ind_{\{t_k > -s\}}\bigg)\bigg]
=\exp\bigg(-\int_0^s F(\upsilon_r(\lambda))\, dr\bigg),
 $$
and
 $$
\P\bigg(\sum_{k \in \ci} \ind_{\{t_k < -s, X_{-t_k}^k > 0\}} = 0 \bigg)
= \exp\bigg(-\int^{\infty}_s F(c(r))\,dr\bigg).
$$

The next equation is obtained by splitting the terms into two parts:
 \begin{eqnarray*}
 \begin{aligned}
&\E_x\bigg[a^{\sum_{j \in J_3}
\ind_{\{X_s^j > 0\}}}\exp\bigg(-\lambda \sum_{j \in J_3} X_s^j\bigg)
\ind_{\{\sum_{j \in J_3} X_s^j > 0\}}\bigg]\\
&\quad = \E_x\bigg[a^{\sum_{j \in J_3}
\ind_{\{X_s^j > 0\}}}\exp\bigg(-\lambda \sum_{j \in J_3} X_s^j\bigg)\bigg]
-\P_x\bigg(\sum_{j \in J_3} X_s^j = 0\bigg).
 \end{aligned}
 \end{eqnarray*}
The first part is calculated as follows:
 \begin{eqnarray*}
 \begin{aligned}
\E_x\bigg[a^{\sum_{j\in J_3}
\ind_{\{X_s^j > 0\}}}\expp{-\lambda \sum_{j \in J_3} X_s^j}\bigg]
&= \exp\bigg(-x \QZ\big[1 - a \ind_{\{X_s > 0\}} \expp{-\lambda X_s} \big]\bigg)\\
& = \exp\bigg(-x \QZ[X_s > 0] + x a \QZ \big[\ind_{\{X_s > 0\}}
\expp{-\lambda X_s}\big]\bigg)\\
&= \exp\bigg(-x \big[(1 - a) c(s) + a \upsilon_s(\lambda)\big]\bigg).
 \end{aligned}
 \end{eqnarray*}
The second part is
 \begin{eqnarray*}
 \begin{aligned}
\P_x \bigg(\sum_{j \in J_3} X_s^j = 0\bigg)
&= \lim_{\lambda \to \infty} \expp{-x \QZ(1 - \expp{-\lambda X_s})}\\
&= \lim_{\lambda\to\infty} \expp{-x \upsilon_s(\lambda)} = \expp{-x c(s)}.
 \end{aligned}
 \end{eqnarray*}
For the last part, we see that
 \beqnn
\QZ[\expp{-\lambda X_s} \ind_{\{X_s > 0\}}]
= \QZ[X_s > 0] -  \QZ[1- \expp{-\lambda X_s}] = c(s)- \upsilon_s(\lambda).
 \eeqnn

Putting all the calculations together we see that
 \begin{eqnarray*}
 \begin{aligned}
\E\big[a^{N^A} \expp{-\lambda Z_0} f(A) \big]
&= \int_0^{\infty} ds\, f(s)
\exp\bigg(-\int_0^s F(\upsilon_r(\lambda))\,dr - \int_s^{\infty} F(c(r))\,dr\bigg)\\
&\qquad\qquad\quad \times \bigg( F(c(s)) - F\big((1 - a) c(s) + a \upsilon_s(\lambda)\big)\bigg).
 \end{aligned}
 \end{eqnarray*}
This finishes the proof.\qed

Using the density of $A$, the following corollary is immediate.
 \begin{cor}\label{c5.1}
For $0 < a < 1,\lambda, t \geq 0$,  we have
 \beqnn
\qquad\E\big[a^{N^A} \expp{-\lambda Z_0} | A = t\big] =
\frac{F(c(t)) - F\big((1 - a)c(t) + a \upsilon_t(\lambda)\big)}
{F(c(t))} \exp\bigg(-\int_0^t F(\upsilon_r(\lambda))\,dr\bigg)
 \eeqnn
and
 \beqlb
\E[ a^{N^A} | A = t] = \frac{F(c(t)) - F\big((1 - a)c(t)\big)}{F(c(t))}
= 1 - \frac{F((1 - a) c(t))}{F(c(t))}.
 \eeqlb
 \end{cor}

The next corollary is direct from Corollary~\ref{c5.1}.
 \begin{cor}\label{c5.2}
We have for $n \geq 1$,
 \beqnn
\P[N^A = n | A = t]=(-1)^{n + 1} \frac{c(t)^n F^{(n)}(c(t))}{n!\ F(c(t))}.
 \eeqnn
Then
$\E[N^A | A = t] = \frac{F^{'}(0+) c(t)}{F(c(t))} \in [0, \infty]$.
In addition if $F^{'}(0+) < \infty$, the function $t\mapsto \EF[N^A | A = t]$ is non-increasing.
 \end{cor}

Notice that if we let $F(t) = c t^{\alpha}, c>0$ and $0 < \alpha < 1$,
then the conditional distribution will
not depend on the CB-process but only on the immigration structure.

\section{The number of ancestors at a fixed time\label{s6}}

In this section we will consider the number of ancestors $M_s$ at time
$-s$ of the current population living at time $0$ and how fast it tends to infinity.
To answer this question we need to introduce the genealogy
of the families which is a richer structure studied in \cite{dl02, dl05, adh12, ad12}.

\subsection{Genealogy of CB-process}

The construction developed by Duquesne and
Le Gall \cite{dl02, dl05} for (sub)-critical CB-process is well known. Results in \cite{dl05}
is restated in the framework of the measured rooted real trees, see \cite{adh12}.
We will follow Section 2 in \cite{ad12}.

\subsubsection{Real tree framework}

A metric space $(\Tau, d)$ is a real tree if the following
two properties hold for every $s, t \in \Tau$,
 \begin{itemize}
\item (unique geodesic) There is a unique isometric map $f_{s, t}$ from $[0,d(s, t)]$
     into $\Tau$ such that $f_{s, t}(0) = s$ and $f_{s, t}(d(s, t)) = t$.
\item (no loop) If $q$ is a continuous injective map from [0,1] into $\Tau$
      such that $q(0) = s$ and $q(1) = t$, we have $q([0, 1]) = f_{s, t}([0, d(s, t)])$.
 \end{itemize}
A rooted real tree is a real tree  $(\Tau, d)$ with a
distinguished vertex $\emptyset$ called the root. Denote such a tree by $(\Tau, d, \emptyset)$.
If $s, t \in \Tau$, we will note $\llbracket s, t \rrbracket$ the range of the
isometric map  $f_{s, t}$ described above.  We also denote
$\llbracket s, t \llbracket = \llbracket s,t \rrbracket \setminus \{t\}$.

If $x \in \Tau$, the degree of $x$, $n(x)$, is the number of connected components
of the set $\Tau \setminus \{x\}$.   The set of leaves is defined as
$\mathrm{Lf}(\Tau) = \{x \in \Tau \backslash \{\emptyset\},\ n(x) = 1\}.$
The skeleton of $\Tau$ is  the set of points in the tree that are not
leaves: $\mathrm{Sk}(\Tau) = \Tau \backslash  \mathrm{Lf}(\Tau)$.

For every $x\in \Tau$, $\lb \emptyset, x \rb$ is
interpreted as the ancestral line of vertex $x$ in the tree.
 If $x, y \in \Tau$, there exists a unique $z \in \Tau$, called the Most Recent
Common Ancestor (MRCA) of $x$ and $y$,  such that
$\lb \emptyset, x \rb \cap \lb \emptyset, y \rb = \lb \emptyset, z \rb$.
Then the root can be seen as the ancestor of all
the population in the tree. We shall call the height of $x$,
 $h(x)$, the distance $d(\emptyset,x)$ to the
root. The function $x\mapsto h(x)$ is continuous on $\Tau$, and we define the height
of $\Tau$ by $H_{max}(\Tau) = \sup_{x\in \Tau} h(x).$

\subsubsection{Measured rooted real trees}

We   will  denote   by   $\T$  the   set   of  the measured rooted real trees
$(\Tau, d, \emptyset, \bm)$  where $(\Tau, d, \emptyset)$  is a locally  compact rooted
real tree and $ \bm$ is a locally finite measure on $\Tau$.  We
may simply write  $\Tau$ in case of  no confusion.

Let $\Tau \in \T$.  For $a \geq 0$, we set
$\Tau(a) = \{x \in \Tau, \, d(\emptyset, x) = a\}$ for the level set at height $a$,
and $\pi_a(\Tau) = \{ x \in \Tau, \ d(\emptyset, x) \le a \}$ for the truncated tree $\Tau$
up to level  $a$. We consider $\pi_a(\Tau)$ with the
root  $\emptyset$,  $d^{\pi_a(\Tau)}$ and $\bm^{\pi_a(\Tau)}$ are the restrictions of $d$ and $\bm$ to
$\pi_a(\Tau)$. Let $(\Tau^{k, \circ}, k \in \ck)$ be the connected components of
$\Tau \setminus  \pi_a(\Tau)$.   Denote by  $\emptyset_k$  the MRCA of all the
vertices  of  $\Tau^{k, \circ}$.   Set $\Tau^k = \Tau^{k, \circ}\cup \{\emptyset _k\}$
which is a real tree rooted at point $\emptyset_k$ with
mass  measure  $\bm^{\Tau^k}$ defined  as  the  restriction  of $\bm^{\Tau}$  to
$\Tau^k$. We will consider the point measure on $\Tau \times \T$:
$$
\cn_a^{\Tau} = \sum_{k \in \ck} \delta_{(\emptyset_k, \Tau^k)}.
$$

\subsubsection{Excursion measure of L\'evy tree}
Recall that $\psi$ is a (sub)-critical branching mechanism.
There exists a $\sigma$-finite measure (or an excursion measure  of Lévy tree)
$\N[d\ct]$ on $\T$, with the
following properties:
\begin{enumerate}
\item[(i)] (Height). $\forall a > 0$,
     $\N[H_{\text{max}}(\ct) > a] = c(a)$.

\item[(ii)] (Mass measure). The mass measure $\bm^\ct$ is supported on
$\mathrm{Lf}(\ct)$, $\N[d\ct]$-a.e.

 \item[(iii)] (Local time). There exists a $\ct$-measure valued process $(\ell^a, a \geq 0)$
càdlàg for the weak topology on finite measure on $\ct$ such that
$\N[d\ct]$-a.e.:
$$
\bm^{\ct}(dx) = \int_0^\infty \ell^a(dx) \, da,
$$
$\ell^0 = 0$, $\inf\{a > 0 ; \ell^a = 0\} = \sup\{a \geq 0 ; \ell^a \neq
0\} = H_{\text{max}}(\ct)$ and for every fixed $a > 0$,
$\N[d\ct]$-a.e.:
 \begin{itemize}
 \item The measure $\ell^a$ is supported on
$\ct(a)$.
 \item We have for every bounded
continuous function $\phi$ on $\ct$:
\begin{align*}
\langle \ell^a, \phi \rangle
& = \lim_{\epsilon \downarrow 0}
\frac{1}{v(\epsilon)} \int \phi(x) \ind_{\{H_{\text{max}}(\ct') \ge
\epsilon\}} \cn_a^{\ct}(dx, d\ct') \\
 & = \lim_{\epsilon \downarrow 0} \frac{1}{v(\epsilon)} \int \phi(x)
\ind_{\{H_{\text{max}}(\ct') \ge \epsilon\}}
\cn_{a-\epsilon}^{\ct}(dx, d\ct').
\end{align*}
 \end{itemize}
 Under $\N$,  the process $(\langle \ell^a, 1 \rangle ,
 a \geq  0)$ is  distributed as  $X$
 under $\QZ$.
\item[(iv)] (Branching property). For every $a>0$, the conditional
distribution of the point measure $\cn_a^{\ct}(dx,d\ct')$ under
$\N[d\ct|H_{\text{max}}(\ct) > a]$, given $\pi_a(\ct)$, is that of a Poisson
point measure on $\ct(a) \times \T$ with intensity
$\ell^a(dx) \N[d\ct']$.
\end{enumerate}

In order to simplifty notations, we will identify $X$ with $(\langle \ell^a, 1 \rangle ,
 a \geq  0)$ as well as $\QZ$ with $\N$.

We give a definition for the number of ancestors.
\begin{defn}
The number of ancestors at time $a$ of the population
living at time $b$ is the number of subtrees above level $a$ which reach level $b > a$:
 \beqnn
R_{a, b}(\Tau) = \sum_{k \in \ck} \ind_{\{H_{max}(\Tau^k) \geq b-a\}}.
 \eeqnn
\end{defn}

\subsection{Genealogy of stationary CBI-process}

We use \eqref{eq:familyd} to construct the genealogy of $Z_0$.

 \begin{itemize}
  \item Conditionally on $N_0$, let
      $\tilde{N}_1(dt,d\Tau) = \sum_{j \in \cj_1} \delta_{(t_j, \Tau^j)}(dt, d\Tau)$
      be a Poisson point measure with intensity $\upsilon(dt) \QZ(d\Tau)$
      with $\upsilon(dt) = \sum_{i \in I} r_i  \delta_{t_i}(dt)$.
  \item Let
      $\tilde{N}_2(dt, d\Tau) = \sum_{j  \in \cj_2}\delta_{(t_j, \Tau^j)}(dt, d\Tau)$
      be a Poisson point measure independent of $N_0, \tilde{N}_1$
      with intensity $\beta dt \QZ(d\Tau)$.
 \end{itemize}
We will write $X^j$ for $\ell^a(\Tau^j)$ for $j \in \cj$. Thus notation \eqref{eq:familyd} is
still consistent with the previous sections. $\sum_{j\in\cj}\delta_{(t_j, \Tau^j)}$
allows to code the genealogy of the family of $Z_0$ in the L\'evy tree sense.

Let $s>0$, we will consider the number of ancestors at time $-s$ of
the current population living at time 0, that is,
$$ M_s=\sum_{j\in\cj}1_{\{t_j<-s\}}R_{-s-t_j, -t_j}(\Tau^j).$$

\subsection{Asymptotic  for the number of ancestors}

First we present the following theorem.
 \begin{thm}\label{t6.1}
The conditional joint distribution of $M_s$ and $Z_0$ is:
for $\eta, \lambda \geq 0, s > 0$,
 $$
\E[\expp{-\eta M_s - \lambda Z_0} | Z_{-s}]
= \expp{-\int_0^s F(\upsilon_r(\lambda))\, dr }
\expp{-Z_{-s}[(1 - \expp{-\eta}) c(s) + \expp{-\eta} \upsilon_s(\lambda)]}.
$$

In particular, conditionally on $Z_{-s}$, $M_s$ is distributed as a
Poisson random variable with parameter $c(s) Z_{-s}$.
 \end{thm}
\pf For any $\eta, \lambda \geq 0$, we have
 \begin{eqnarray*}
 \begin{aligned}
&\E[\expp{-bZ_{-s} - \eta M_s - \lambda Z_0}]\\
&= \E\left[\exp\bigg(-\lambda \sum_{j \in \cj} \ind_{\{-s < t_j \leq 0\}}
X_{-t_j}^j\bigg)\right]\E\left[\exp\bigg(-bZ_{-s} - \eta M_s
-\lambda \sum_{j \in \cj} \ind_{\{ t_j \leq -s\}} X _{-t_j}^j\bigg)\right]\\
&= \exp\bigg(-\int_0^s F\big(\upsilon_r(\lambda)\big) dr \bigg)\\
&\qquad\qquad\qquad\qquad\times \E\bigg[\exp\bigg(-\sum_{j \in \cj}
\ind_{\{ t_j \leq -s\}}(b\ X_{-s - t_j}^j + \eta R_{-s - t_j, -t_j}(\Tau^j)
+ \lambda X _{-t_j}^j)\bigg)\bigg]\\
&= \exp\bigg(-\int_0^s F(\upsilon_r(\lambda))\,dr \bigg)\exp\bigg(-\int_0^{\infty} da\,
F \big[\QZ\big(1 - \expp{-bX_a - \eta R_{a, a + s}(\Tau) - \lambda X_{a + s}}\big)\big]\bigg).
 \end{aligned}
 \end{eqnarray*}
where in the first equality we use the fact that the Poisson point
measures over disjoint sets are independent, in the second one we use
Lemma \ref{l3.1} and  we use an immediately generalization of Lemma \ref{l3.1}
to genealogies in the third equality.
Using branching property  we have
 \begin{eqnarray*}
 \begin{aligned}
\QZ\left[1 - \expp{-b X_a - \eta R_{a, a + s}(\Tau) - \lambda X_{a + s}}\right]
&= \QZ\left[1 - \expp{-b X_a-\sum_{k \in \ck}(\eta \ind_{\{H_{max}(\Tau^k) \geq s\}} - \lambda X(\Tau^k)_s)}\right]\\
&= \QZ\left[1 - \expp{-X_a \{b + \QZ[1 - \exp(-\eta \ind_{\{H_{max}(\Tau) \geq s\}}-\lambda X_s)]\}}\right].
 \end{aligned}
 \end{eqnarray*}

Since  $X_s = 0$ on $\{H_{max}(\Tau)< s\}$, we have $1 - \expp{-\eta \ind_{\{H_{max}(\Tau) \geq s\}} - \lambda X_s}
= (1 - \expp{-\eta})\ind_{\{\zeta \geq s\}} + \expp{-\eta}(1 - \expp{-\lambda X_s})$.
Then we deduce that
$$
\QZ\left[1 - \expp{-b X_a - \eta R_{a, a + s}(\Tau) - \lambda X_{a + s}}\right]
= \QZ\big[1 - \expp{-\lambda^{'} X_a}\big]
= \upsilon_a(\lambda^{'})
$$
with $\lambda^{'} = b + (1 - \expp{-\eta}) c(s) + \expp{-\eta}\upsilon_s(\lambda)$.
Then we get the result. \qed

Intuitively $M_s$ counts the number of excursions of the height process
at time $-s$ above level $s$. Similar result as that of Duquesne and Le
Gall \cite{dl05} can be deduced here.

 \begin{cor}\label{c6.2}
The following convergence holds:
$$
\lim_{s \to 0} \frac{M_s}{c(s)} = Z_0, \quad \mbox{a.s.}
$$
 \end{cor}

\section{The MRCA age process and the zero set of the CBI-process}

In this section we will deal with the MRCA age process $(A_t, t \in \R)$
and the zero set of the CBI-process $\cz = \overline {\{t \in \R, Z_t = 0\}}$
by using the sample path decomposition of the CBI-process.

For the clan decomposition of the CBI-process shown in Section 3,
we have used the Poisson point process
$N_3(dt, dX) = \sum_{i \in \ci} \delta_{(t_i, X^i)}(dt, dX)$
with intensity $dt \mu(dX)$, where
$\mu(dX) = \beta \QZ(dX) + \int_{(0,\infty)} n(dx) \P_x(dX).$
The corresponding Poisson point process which is in charge of the birth time
and duration time is $N_4(dt, d\zeta) = \sum_{i \in \ci} \delta_{(t_i, \zeta^i)}(dt, d\zeta)$
with intensity $dt \mu(\zeta \in \cdot)$, where $\mu(\zeta > t) = F(c(t))$.

\subsection{The MRCA age process}

Define the left leaning wedge with apex at $(t,r)$ by
$$
\Delta(t,r):= \{(u,v), u < t\, \mbox{and} \, u+v > r\},
$$
which is the set of points that give birth before $t$ and is still alive at time $r$.
We can thus define the MRCA age process $(A_t, t \in \R)$ as follows:
$$
A_t:=t- \inf\{s: \exists \zeta > 0,\, \mbox{such\, that}\,
(s,\zeta)\in \{(t_i, \zeta_i): i \in \ci\} \cap \Delta(t,0) \}.
$$
The strong Markov property of the Poisson point processes $N_4$ implies
that $(A_t, t \in \R)$ is a time homogeneous Markov process.
The transition probabilities of the MRCA age process can be proved
 along the same lines as that of part (a) of Theorem 1.1 in Evans and Ralph \cite{er10}.
\begin{thm}
The transition probabilities of the Markov process $(A_t, t \in \R)$ are seperated into two parts:
\begin{itemize}
\item for $0< y < x + t$,  the continuous part
$$
\P(A_{s+t} \in dy | A_s = x) = \left(1 - \frac{F(c(x + t))}{F(c(x))}\right )
\expp{-\int_y^{x + t} F(c(u))\,du} F(c(y)) \, dy,
$$
\item otherwise, a single atom
$\P(A_{s + t} = x + t | A_s = x) = \frac{F(c(x + t))}{F(c(x))}.$
\end{itemize}
\end{thm}

\subsection{The zero set of the CBI-process}

From the clan decomposition,  we have

\begin{itemize}
  \item either $Z_t \neq 0$, for $t \in [t_i, t_i + \zeta_i), i \in I$;
  \item or $Z_t \neq 0$, for $t \in (t_i, t_i + \zeta_i),i \in \cj_2$.
\end{itemize}
Then we have
$$
\{t, Z_t = 0\} = \R \setminus \left[\bigcup_{i \in I}[t_i, t_i + \zeta_i) \cup
\bigcup_{i \in \cj_2} (t_i, t_i + \zeta_i)\right].
$$
Then we can derive that the zero set of the CBI-process is a random renewal set with the following equation:
$$
\cz = \R \setminus \bigcup_{i \in \ci}(t_i, t_i + \zeta_i).
$$
To see why this is always true, we only need to focus on those points that belong to
$\R \setminus \bigcup_{i \in \ci}(t_i, t_i + \zeta_i)$ but not to $\{t, Z_t = 0\}$.
These points actually are the left accumulation points of $\{t, Z_t = 0\}$. Indeed suppose that this is not true.
Let $\{t_i\}$ be such point. Then for any $t_i$, there exists $\epsilon_i > 0$ such that
$[t_i - \epsilon_i, t_i)\cap \{t, Z_t = 0\}=\emptyset.$  This is impossible since after taking the closure
we obtain that $[t_i - \epsilon_i, t_i]\cap \cz = \emptyset;$ while $t_i \in \cz$.

We derive the following theorem.
\begin{thm}
\begin{enumerate}
  \item $\cz = \emptyset$ if and only if $\int_0^1 \exp{(\int_t^\infty F(c(u)) du)}\, dt = \infty$.
  \item If $\int_0^1 \exp{(\int_t^\infty F(c(u)) du)}\, dt < \infty$,
   the radom set $\cz$ has a positive Lebesgue measure a.s. if and only if
   $\int_0^\infty \frac{F(t)}{\psi(t)}\, dt < \infty$;  the random set $\cz$ is the union of the
   closed intervals of positive lengths  if and only if $\beta = 0$ and $n(dx) < \infty$.
\end{enumerate}
\end{thm}
\pf We have known that $\mu(\zeta > t) = F(c(t))$.
We can derive directly from Corollary 5 in \cite{ffs85} that the first assertion holds.

With a slightly modification of Proposition 1 in \cite{ffs85} or Propositon 1.22 in \cite{ffm85},
we can obtain that $\cz$ has a positive Lesbesgue measure if and only if
$$
\int_0^\infty t\, \mu(\zeta \in dt)=\int_0^\infty \mu(\zeta > t) dt = \int_0^\infty  F(c(t)) dt < \infty.
$$
Letting $r = c(t)$ in the above equation  yields $\int_0^\infty \frac{F(t)}{\psi(t)}\, dt < \infty$.
In order to derive the last statement, we use Corollary 2 in \cite{ffs85},
which requires $\mu(\zeta \in dt)$ to be finite, i.e. $F(c(t)) $ is finite, as $t \to 0+$.
Since we have $c(t) \to \infty$ as   $t \to 0+$,   we need $\beta = 0$ and $n((0, \infty)) < \infty$.\qed

Easy calculation gives a simple example that $\psi(u)=2 \beta \theta u + \beta u^2$ and $F(u)= 2 \beta u$ satisfying
condition $(1)$, and $\cz = \emptyset$. Another example is given in \cite{fh13} in stable case.

{\bf Acknowledgement.} I would like to express my sincere thanks to my advisor Professor Zenghu Li for his
persistent encouragement, suggestions. I also want to thank Professor J.F. Delmas
for his careful check of this work and the suggestions.

\newcommand{\sortnoop}[1]{}


\begin{thebibliography}{1}

\bibitem{ad12}
R. Abraham and J.-F. Delmas.
\newblock The forest associated with the record process on a L\'evy tree.
\newblock {\em arXiv:1204.2357}, 2012.
\bibitem{adh12}
R. Abraham, J.-F. Delmas and P. Hoscheit.
\newblock Exit times for an increasing L\'evy tree-valued process.
\newblock  {\em arxiv: 1202.5463 }, 2012.

\bibitem{a91}
D. Aldous.
\newblock  The continuum random tree I.
\newblock  {\em Ann.  Probab.}, 19(1): 1-28, 1991.

\bibitem{a93}
D. Aldous.
\newblock The continuum random tree III.
\newblock {\em Ann. Probab.}, 21(1): 248-289, 1993.

\bibitem{Ali85}
S. A. Aliev.
\newblock  A limit theorem for the  Galton--Watson branching processes with immigration.
\newblock  {\em Ukrainian Math. J.}, 37: 535--438, 1985.

\bibitem{AlS82}
S. A. Aliev and V.M. Shchurenkov.
\newblock    Transitional phenomena and the convergence of Galton--Watson processes to Ji\v{r}ina processes.
\newblock {\em Theory Probab. Appl.}, 27: 472--485, 1982.


\bibitem{bbs11}
J. Berestycki, N. Berestycki and V. Limic.
\newblock A small time coupling between $\Lambda$-coalescent and branching processes.
\newblock {\em arXiv: 1101.1875}, 2011.


\bibitem{cd10}
Y.-T. Chen and J.-F. Delmas.
\newblock Smaller population size at the MRCA time for stationary branching processes.
\newblock {\em Ann. Probab. To appear (arXiv: 1009.0814),}  2012.


\bibitem{dl02}
T. Duquesne and J.-F. Le Gall.
\newblock {\it Random trees, L\'evy processes and spatial branching processes.}
\newblock Volume 281. Ast\'erisque, 2002.



\bibitem{dl05}
T. Duquesne and J.-F. Le Gall.
\newblock Probabilistic and fractal aspects of L\'evy trees.
\newblock {\em Probab. Theory  Relat. Fields},  131(4): 553-603, 2005.


\bibitem{er10}
S. N. Evans and P. L. Ralph.
\newblock Dynamics of the time to the most recent common ancestor in a large branching population.
\newblock {\em Ann. Appl. Probab.},  20(1): 1-25, 2010.

\bibitem{fb12}
C. Foucart and G. U. Bravo.
\newblock Local extinction in continuous state branching processes with immigration.
\newblock  {\em arXiv: 1211.3699}, 2012.

\bibitem{fh13}
C. Foucart and O. H\'enard.
\newblock Stable continuous-state branching processes with immigration and Beta-Fleming-Viot
processes with immigration.
\newblock  {\em Electron. J. Probab.}, 23: 1-21, 2013.

\bibitem{ffm85}
P. J. Fitzsimmons, B. Fristedt and B. Maisonneuve.
\newblock Intersections and limits of regenerative sets.
\newblock {\em Z.Wahrsch. Verw. Gebiete}, 70: 157-173, 1985.


\bibitem{ffs85}
P. J. Fitzsimmons, B. Fristedt and L. A. Shepp
\newblock The set of real numbers left uncovered by random covering intervals.
\newblock {\em Z.Wahrsch. Verw. Gebiete},  70: 175-189, 1985.

\bibitem{ft88}
P. J. Fitzsimmons and M. Taksar.
\newblock Stationary regenerative sets and subordinators.
\newblock {\em Ann. Probab.}, 16(3): 1299-1305, 1988.

\bibitem{J58}
 M. Jirina.
\newblock Stochastic branching processes with continuous state space.
\newblock {\em Czech. Math. J.}, 83(8): 292-312, 1958.

\bibitem{kw71}
K. Kawazu and S. Watanabe.
\newblock Branching processes with immigration and related limit theorems.
\newblock {\em Theory Probab. Appl.},  83(8): 36-54, 1971.


\bibitem{k82a}
J. F. C. Kingman.
\newblock On the genealogy of large populations.
\newblock {\em J. Appl. Probab.},  19: 27-43, 1982.


\bibitem{k82b}
J. F. C. Kingman.
\newblock The coalescent.
\newblock {\em  Stochastic Process. Appl.},  13(3): 235-248, 1982.



\bibitem{l07}
A. Lambert.
\newblock Coalescence times for the branching process.
\newblock {\em Adv. Appl. Probab.},  35(4): 1071-1089, 2007.


\bibitem{l67}
J. Lamperti.
\newblock The limit of a sequence of branching processes.
\newblock {\em Probab. Theory  Relat. Fields},  7(4): 271-288, 1967.


\bibitem{Li06}
Z. Li.
\newblock A limit theorem for discrete Galton--Watson branching processes with immigration.
\newblock {\em J. Appl. Probab.}, 43: 289--295, 2006.

\bibitem{l11}
Z. Li.
\newblock {\it Measure-valued branching Markov processes}.
\newblock Springer, Berlin, 2011.


\bibitem{l12}
Z. Li.
\newblock {\it Continuous-state branching processes.}
\newblock Lecture Notes, Beijing Normal University. {\em arXiv:1202.3223}, 2012.

\bibitem{p99}
J. Pitman.
\newblock Coalescents with Multiple collisions.
\newblock {\em Ann. Probab.}, 27: 1870-1902, 1999.

\bibitem{s99}
S. Sagitov.
\newblock The general coalescent with asynchronous mergers of ancestral lines.
\newblock {\em J. Appl. Probab.},  36: 1116-1125, 1999.

\bibitem{t80}
M. I. Taksar
\newblock regenerative sets on real line.
\newblock{\em Seminaire de Probabilites XIV in Lecture Notes in Mathematics.},  784: pp. 437-474. Springer, New York, 1980.


\end{thebibliography}
 \end{document}